\documentclass[12pt]{article}
\usepackage{enumerate}
\usepackage{latexsym}
\usepackage{amsfonts}
\usepackage[only,ninrm,elvrm,twlrm,sixrm,egtrm,tenrm]{rawfonts}
\usepackage{indentfirst}
\usepackage{amsmath}
\usepackage{graphicx,psfrag}
\usepackage{graphics}
\usepackage{makeidx}
\newtheorem{thm}{Theorem}[section]

\newtheorem{lem}[thm]{Lemma}

\newtheorem{prop}[thm]{Proposition}

\def\qed{\hfill \nopagebreak\rule{5pt}{8pt}}
\def\pf{\noindent {\it Proof.} }
\title{\bf Bipartite rainbow numbers of matchings
\footnote {Supported by NSFC and the ``973" project.}}

\author{
\small Xueliang Li and Jianhua Tu\\
\small Center for Combinatorics and LPMC, Nankai University\\
\small Tianjin 300071, P.R. China. Email: lxl@nankai.edu.cn\\
\small Zemin  Jin\\
\small Department of Mathematics, Zhejiang Normal University\\
\small Jinhua, Zhejiang 321004, P.R. China. Email: zeminjin@zjnu.cn
}
\date{}

\begin{document}

\maketitle
\begin{abstract}
Given two graphs $G$ and $H$, let $f(G,H)$ denote the maximum number
$c$ for which there is a way to color the edges of $G$ with $c$
colors such that every subgraph $H$ of $G$ has at least two edges of
the same color. Equivalently, any edge-coloring of $G$ with at least
$rb(G,H)=f(G,H)+1$ colors contains a rainbow copy of $H$, where a
rainbow subgraph of an edge-colored graph is such that no two edges
of it have the same color. The number $rb(G,H)$ is called the {\it
rainbow number of $H$ with respect to $G$}, and simply called the
{\it bipartite rainbow number of $H$} if $G$ is the complete
bipartite graph $K_{m,n}$. Erd\H{o}s, Simonovits and S\'{o}s showed
that $rb(K_n,K_3)=n$. In 2004, Schiermeyer determined the rainbow
numbers $rb(K_n,K_k)$ for all $n\geq k\geq 4$, and the rainbow
numbers $rb(K_n,kK_2)$ for all $k\geq 2$ and $n\geq 3k+3$. In this
paper we will determine the rainbow numbers $rb(K_{m,n},kK_2)$ for
all $k\geq 1$.\\
[0.1in]{\bf Keywords:} edge coloring, rainbow subgraph, rainbow number.\\
\end{abstract}

\section{Introduction}

In this paper we consider undirected, finite and simple graphs only,
and use standard notations in graph theory (see \cite{bon}). Given
two graphs $G$ and $H$, if $G$ is edge-colored and a subgraph $H$ of
$G$ contains no two edges of the same color, then $H$ is called a
{\it totally multicolored (TMC)} or {\it rainbow subgraph} of $G$
and we say that $G$ contains a TMC or rainbow $H$. Let $f(G,H)$
denote the maximum number of colors in an edge-coloring of the graph
$G$ with no TMC $H$. We now define $rb(G,H)$ as the minimum number
of colors such that any edge-coloring of $G$ with at least
$rb(G,H)=f(G,H)+1$ colors contains a TMC or rainbow subgraph
isomorphic to $H$. The number $rb(G,H)$ is called the {\it rainbow
number of $H$ with respect to $G$}. If $G$ is the complete bipartite
graph $K_{m,n}$, $rb(G,H)$ is simply called the {\it bipartite
rainbow number of $H$}.

When $G=K_n$, $f(G,H)$ is called the {\it anti-Ramsey number} of
$H$. Anti-Ramsey numbers were introduced by Erd\H{o}s, Simonovits
and S\'{o}s in the 1970s. Let $P_k$ and $C_k$ denote the path and
the cycle with $k$ vertices, respectively. Simonovits and S\'{o}s
\cite{sos} determined $f(K_n,P_k)$ for large enough $n$. Erd\H{o}s
et al. \cite{erd} conjectured that for every fixed $k\geq 3$,
$f(K_n,C_k)=n(\frac{k-2}{2}+\frac{1}{k-1})+O(1)$, and proved it for
$k=3$ by showing that $f(K_n,C_3)=n-1$. Alon \cite{alon} showed that
$f(K_n,C_4)=\lfloor\frac{4n}{3}\rfloor-1$, and the conjecture is
thus proved for $k=4$. Recently, the conjecture is proved for all
$k\geq 3$ by Montellano-Ballesteros and Neumann-Lara \cite{mont}.
Axenovich, Jiang and K\"{u}ndgen \cite{tao} determined
$f(K_{m,n},C_{2k})$ for all $k\geq 2$.

In 2004, Schiermeyer \cite{Ingo} determined the rainbow numbers
$rb(K_n,K_k)$ for all $n\geq k\geq 4$, and the rainbow numbers
$rb(K_n,kK_2)$ for all $k\geq 2$ and $n\geq 3k+3$, where $H=kK_2$ is
a matching $M$ of size $k$. The main focus of this paper is to
consider the analogous problem for matchings when the host graph $G$
is a complete bipartite graph $K_{m,n}$ (say $m\geq n$). For all
positive integers $m\geq n$ and $k\geq 1$, we determine the exact
values of $rb(K_{m,n},kK_2)$.

\section{Main results}

Let $M$ be a matching in a given graph $G$, then the subgraph of $G$
induced by $M$, denoted by $\langle M\rangle_G$ or $\langle
M\rangle$, is the subgraph of $G$ whose edge set is $M$ and whose
vertex set consists of the vertices incident with some edges in $M$.
A vertex of $G$ is said to be \emph{saturated} by $M$ if it is
incident with an edge of $M$; otherwise, it is said to be
\emph{unsaturated}. If every vertex of a vertex subset $U$ of $G$ is
saturated, then we say that $U$ is saturated by $M$. A matching with
maximum cardinality is called a maximum matching.

In a given graph $G$, $N_G(U)$ denotes the set of vertices of $G$
adjacent to the vertex set $U$. If $R,T\in V(G)$, we denote $E(R,T)$
or $E_G(R,T)$ as the set of all edges having a vertex from both $R$
and $T$. Let $G(m,n)$ denote a bipartite graph with bipartition
$A\cup B$, and $|A|=m$ and $|B|=n$, without loss of generality, in
the following we always assume $m\geq n$.

Let $ext(m,n,H)$ denote the maximum number of edges that a bipartite
graph $G(m,n)$ can have with no subgraph isomorphic to $H$. The
bipartite graphs attaining the maximum for given $m$ and $n$ are
called extremal graphs.

We now determine the value $ext(m,n,H)$ for $H=kK_2$.
\begin{lem}\label{ore}
Let $G(m,n)$ be a bipartite graph with bipartition $A\cup B$, and
$M$ a maximum matching in $G$. Then the size of $M$ is $m-d$, where
$$d=\max\{|S|-|N_G(S)|: S\subseteq A\}.$$
\end{lem}

\begin{thm}\label{ext}
$$ext(m,n,kK_2)=m(k-1)\text{ for all }1\leq k\leq n,$$
that is, for any given bipartite graph $G(m,n)$, if
$|E(G(m,n))|>m(k-1)$, then $kK_2\subset G$. Moreover, $K_{m,(k-1)}$
is the unique such extremal graph.
\end{thm}
\pf Suppose that $G$ contains no $kK_2$. Let $M$ be a maximum
matching of $G$ and the size of $M$ is $k-i$, where $i\geq 1$. By
Lemma \ref{ore}, there exists a subset $S\subset A$ such that
$|S|-|N_G(S)|=m-k+i$. Thus
$$|E(G)|\leq|S||N_G(S)|+n(m-|S|)=(|N_G(S)|+m-k+i)|N_G(S)|+n(k-i-|N_G(S)|).$$

Since $0\leq|N_G(S)|\leq k-i\leq k-1$, we obtain
$$|E(G)|\leq \max\{m(k-1),n(k-1)\}\leq m(k-1),$$
where the equality is possible only if $i=1$ and $G\cong K_{m,k-1}$.
So, $K_{m,k-1}$ is the unique such extremal graph.\qed

For $k=1$, it is clear that $rb(K_{m,n},K_2)$=1. Now we determine
the value of $rb(K_{m,n},2K_2)$ (for $k=2$).

\begin{thm}
$$rb(K_{2,2},2K_2)=3,$$
and
$$rb(K_{m,n},2K_2)=2 \text{ for all }m\geq3 \text{ and } n\geq2.$$
\end{thm}
\pf It is obvious that $rb(K_{2,2},2K_2)\leq 3$. Let
$\{a_1,a_2\}\cup\{b_1,b_2\}$ be the two parts of $K_{2,2}$. If
$K_{2,2}$ is edge-colored with 2 colors such that
$c(a_1b_1)=c(a_2b_2)=1$ and $c(a_1b_2)=c(a_2b_1)=2$, then $K_{2,2}$
contains no TMC $2K_2$. So, $rb(K_{2,2},2K_2)=3$.

For $m\geq 3$ and $n\geq 2$, let the edges of $G=K_{m,n}$ be colored
with at least 2 colors. We suppose that the two parts of $K_{m,n}$
are $A$ and $B$ with $|A|=m$ and $|B|=n$. Suppose that $K_{m,n}$
contains no TMC $2K_2$. Let $e_1=a_1b_1$, $a_1\in A$, $b_1\in B$, be
an edge with $c(e_1)=1$, and $R=V(K_{m,n})-\{a_1,b_1\}$. Then
$c(e)=1$ for all edges $e\in E(G[R])$. Moreover, $c(e)=1$ for all
edges $e\in E(b_1,R)$, since $m\geq 3$. Thus $c(e)=1$ for all edges
$e\in E(a_1,R)$. But then $K_{m,n}$ is monochromatic, a
contradiction. So, $rb(K_{m,n},2K_2)=2$ for all $m\geq 3$ and $n\geq 2$.\qed\\

The next proposition provides a lower and upper bound for
$rb(K_{m,n},kK_2)$.

\begin{prop}
$ext(m,n,(k-1)K_2)+2\leq rb(K_{m,n},kK_2)\leq ext(m,n,kK_2)+1$.
\end{prop}
\pf The upper bound is obvious. For the lower bound, an extremal
coloring of $K_{m,n}$ can be obtained from an extremal graph
$K_{m,k-2}$ for $ext(m,n,(k-1)K_2)$ by coloring the edges of
$K_{m,k-2}$ differently and the edges of $\overline{K_{m,k-2}}$ by
one extra color. So, the coloring does not contain a TMC $kK_2$.\qed\\

We will now show that the lower bound can be achieved for all $m\geq
n\geq k\geq 3$, and thus obtain the exact value of
$rb(K_{m,n},kK_2)$ for all $k\geq 3$.

\begin{thm}
$rb(K_{m,n},kK_2)=ext(m,n,(k-1)K_2)+2=m(k-2)+2$ for all $m\geq n\geq
k\geq 3$.
\end{thm}
\pf For $m\geq n\geq k\geq 3$, let the edges of $K_{m,n}$ be colored
with $m(k-2)+2$ colors. Suppose that $K_{m,n}$ contains no TMC
$kK_2$. Since $m(k-2)+2=ext(m,n,(k-1)K_2)+2$, there is a TMC
$(k-1)K_2$ in the coloring of $K_{m,n}$. Now let $G\subset K_{m,n}$
be a TMC spanning subgraph of size $m(k-2)+2$ containing a
$(k-1)K_2$. We suppose that the two parts of the bipartite graph $G$
are $A$ and $B$ with $|A|=m$ and $|B|=n$. By Lemma \ref{ore}, there
exists a subset $S$ of $A$ such that $|S|-|N_G(S)|=m-k+1$, $0\leq
|N_G(S)|\leq k-1$. First, we prove the following two assertions.\\

\noindent{\textbf{Claim 1}}: If one component of $G$ consists of a
$K_{m,k-2}$ and two adjacent pendant edges and the other components
are isolated vertices (see Figure 1), then $K_{m,n}$ contains a TMC
$kK_2$.
\begin{figure}
\begin{center}
\includegraphics[width=7cm]{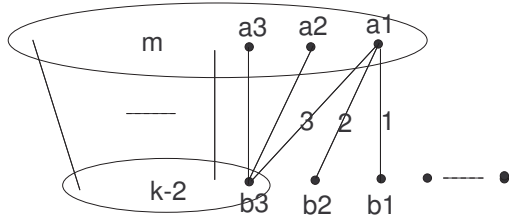}
\caption{The special graph $SG_1$.}
\end{center}
\end{figure}

Denote $SG_1$ as the special graph $G$ and $Q$ as the set of
isolated vertices of $G$. The proof of the claim is given by
distinguishing the following two cases:

Case I. $m\geq 4$.

Without loss of generality, we suppose that $c(a_1b_1)=1$,
$c(a_1b_2)=2$ and $c(a_1b_3)=3$. We will show that $c(a_2b_2)=1$. If
$c(a_2b_2)=2$ or 3, it is obvious that there is a TMC $kK_2$ in
$K_{m,n}$. Otherwise, we suppose $c(a_2b_2)=4$. In $G_1=G-\{Q\cup
a_1\cup a_2\cup b_1\cup b_2\}$, the number of edges whose colors are
not 4 is at least $(m-2)(k-2)-1$. Since $m\geq 4$, we have
$(m-2)(k-2)-1\geq ext(m-2,k-2,(k-2)K_2)+1=(m-2)(k-3)+1$. Thus we can
obtain a TMC $H=(k-2)K_2$ in $G_1$, and hence there is a TMC
$kK_2=H\cup \{a_1b_1,a_2b_2\}$ in $K_{m,n}$. So $c(a_2b_2)$ must be
1. By the same token, $c(a_3b_1)$ must be 2. Now we can obtain a TMC
$H'=(k-3)K_2$ in $G_2=G-\{Q\cup a_1\cup a_2\cup a_3\cup b_1\cup
b_2\cup b_3\}$, and hence there is a TMC $kK_2=H'\cup
\{a_1b_3,a_2b_2,a_3b_1\}$ in $K_{m,n}$.

Case II. $m=n=k=3$.

Without loss of generality, we suppose that $c(a_1b_1)=1$,
$c(a_1b_2)=2$, $c(a_1b_3)=3$, $c(a_2b_3)=4$ and $c(a_3b_3)=5$.
Suppose $K_{m,n}$ contains no $3K_2$. Hence $c(a_2b_2)\in\{1,5\}$
and $c(a_3b_1)\in\{4,2\}\cap \{3,c(a_2b_2)\}=\emptyset$, a contradiction.\\

\noindent{\textbf{Claim 2}}: If one component of $G$ consists of a
$K_{m-1,m-1}$ and a pendant edge (say $pv$ with $d_G(v)=1$) and the
other component is an isolated vertex (say $u$), then $K_{m,n}$
contains a TMC $mK_2$.

Denote $SG_2$ as the special graph $G$. The proof of the claim is
given as follows:

Without loss of generality, we suppose $c(uv)=1$. Then in
$G_3=G-u-v$, the number of edges whose colors are not 1 is at least
$(m-1)(m-1)-1$. Since $m\geq 3$, we have $(m-1)(m-1)-1\geq
ext(m-1,m-1,(m-1)K_2)+1=(m-1)(m-2)+1$. Thus we can obtain a TMC
$H=(m-1)K_2$ in $G_3$, and hence there is a TMC $mK_2=H\cup uv$ in
$K_{m,n}$.\\

We define $R_1=S$, $R_2=A-S$, $T_1=N_G(S)$, $T_2=B-N_G(S)$. Let $M$
be a maximum matching of $G$, then $T_1$ and $R_2$ are saturated by
$M$. There exists an edge $e=z_1z_2\in E_{K_{m,n}}(R_1,T_2)$,
$z_1\in R_1$, $z_2\in T_2$, $z_1\notin\langle M\rangle,z_2\notin
\langle M\rangle$, without loss of generality, we suppose $c(e)=1$.
See Figure 2. So, there exists an edge $e_1\in M$ such that
$c(e_1)=1$. Now we distinguish three cases to finish the proof of
the theorem.

\begin{figure}
\begin{center}
\includegraphics[width=7cm]{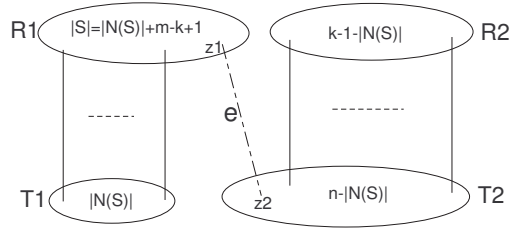}
\caption{$T_1$ and $R_2$ are saturated by $M$.}
\end{center}
\end{figure}

Case 1. $|N_G(S)|=k-1$.

In this case, $R_1=A$, there is no $(k-1)K_2$ in $G'=G-(T_2\cup
z_1)-e_1$. By Theorem \ref{ext}, $|E(G')|\leq (m-1)(k-2)$. Thus,
$$|E(G)|=1+|E(G')|+|E_G(z_1,T_1)|\leq1+(m-1)(k-2)+(k-1)\leq m(k-2)+2.$$
If $m>k$, since $|E(G)|=m(k-2)+2$, then $G'=K_{m-1,k-2}$ and
$|E_G(z_1,T_1)|=k-1$. It is easy to check that $(G-e_1+e)\cong
SG_1$, and by the proof of {\bf Claim 1} we can obtain a TMC $kK_2$
in $K_{m,n}$. If $m=k$, again since $|E(G)|=m(k-2)+2$, it is easy to
check that $(G-e_1+e)\cong SG_1$ or $G\cong SG_2$, and by {\bf Claim
1} and {\bf Claim 2} we can obtain a TMC $kK_2$ in $K_{m,n}$.

Case 2. $|N_G(S)|=0$.

In this case, $G'=G-(R_1\cup z_2)-e_1$ and there is no $(k-1)K_2$ in
$G'$. Similarly,
$$|E(G)|=1+|E(G')|+|E_G(z_2,R_2)|\leq1+(n-1)(k-2)+(k-1)\leq n(k-2)+2.$$
If $m>n$, this contradicts that $G$ has $m(k-2)+2$ edges; if $m=n$,
by Case 1 we can obtain a TMC $kK_2$ in $K_{m,n}$.

Case 3. $1\leq |N_G(S)|\leq k-2$.

\indent{Subcase 3.1}. $e_1\in E_G(R_2,T_2)$.

In this case, there is no TMC $(k-1-|N_G(S)|)K_2$ in $G'=G[R_2\cup
T_2]-z_2-e_1$. Thus,
\begin{eqnarray*}
|E(G)|&=&1+|E(G')|+|E_G(z_2,R_2)|+|E_G(T_1,A)|\\
&\leq&1+(k-2-|N_G(S)|)(n-|N_G(S)|-1)\\
&&+(k-1-|N_G(S)|)+m|N_G(S)|\\
&\leq&\max\{3+n(k-2)+(m-n)-(k-2),m(k-2)+2\}\\
&\leq&m(k-2)+2.
\end{eqnarray*}
Since $|E(G)|=m(k-2)+2$, it is easy to check that $G'$ is an empty
graph and $G\cong SG_1$, and hence there is a TMC $kK_2$ in
$K_{m,n}$.

\indent{Subcase 3.2}. $e_1\in E_G(R_1,T_1)$.

In this case, $G'=G[R_1\cup T_1]-z_1-e_1$ and there is no TMC
$|N_G(S)|K_2$ in $G'$. Similarly,
\begin{eqnarray*}
|E(G)|&=&1+|E(G')|+|E_G(z_1,T_1)|+|E_G(R_2,B)|\\
&\leq&1+(|N_G(S)|-1)(|N_G(S)|+m-k)+|N_G(S)|+n(k-1-|N_G(S)|)\\
&\leq&\max\{3+m(k-2)+(n-m)-(k-2),n(k-2)+2\}\\
&\leq&m(k-2)+2.
\end{eqnarray*}
If $m>n$, then $|E(G)|<m(k-2)+2$, a contradiction. Otherwise,
$|E(G)|=m(k-2)+2$ only if $|N_G(S)|=1$ and $G'$ is an empty graph
and $G\cong SG_1$, hence there is a TMC $kK_2$ in $K_{m,n}$.

The proof is now complete.\qed


\begin{thebibliography}{99}

\bibitem{alon} N. Alon, On a conjecture of Erd\H{o}s, Simonovits and
S\'{o}s concerning anti-Ramsey theorems, {\it J. Graph Theory} {\bf
7}(1983), 91-94.

\bibitem{tao} M. Axenovich, T. Jiang, and A. K\"{u}ndgen,
Bipartite anti-Ramsey numbers of cycles, {\it J. Graph Theory} {\bf
47}(2004), 9-28.

\bibitem{bon} J.A. Bondy and U.S.R. Murty, ``Graph Theory with
Applications", Macmillan, London; Elsevier, New York, 1976.

\bibitem{erd} P. Erd\H{o}s, M. Simonovits and V.T. S\'{o}s,
Anti-Ramsey theorems, in: A. Hajnal, R. Rado, V.T. S\'{o}s (Eds),
Infinite and Finite Sets, Vol. II, Colloq. Math. Soc. J\'{a}nos
Bolyai {\bf 10}(1975), 633-643.

\bibitem{mont} J.J. Montellano-Ballesteros and V. Neumann-Lara,
An anti-Ramsey theorem on cycles, {\it Graphs and Combinatorics}
{\bf 21}(2005), 343¨C354.

\bibitem{Ingo} I. Schiermeyer, Rainbow numbers for matchings and
complete graphs, {\it Discrete Math.} {\bf 286}(2004), 157-162.

\bibitem{sos} M. Simonovits and V.T. S\'{o}s, On restricted
colourings of $K_n$, {\it Combinatorica} {\bf 4}(1984), 101-110.

\end{thebibliography}
\end{document}